\documentclass[12pt]{article}
\usepackage[latin1]{inputenc}
\usepackage{hyperref}
\usepackage{amsmath}
\usepackage{amsfonts}
\usepackage{amssymb}
\usepackage{graphicx}
\usepackage[mathscr]{eucal} % usa the other fuente (euler) in the equations
\usepackage{subfigure}
\usepackage{enumerate}
\usepackage{amsthm} 
\usepackage{xcolor}
\usepackage{float}
\usepackage[boxed,vlined]{algorithm2e}

\graphicspath{ {../Figures/} }

\theoremstyle{plain}

\newtheorem{theorem}{Theorem}
\newtheorem{proposition}[theorem]{Proposition}
\newtheorem{definition}[theorem]{Definition}

\newcommand{\eqsepv}{\; , \enspace}       
\newcommand{\eqfinv}{\; ,}                
\newcommand{\eqfinp}{\; .}

\newcommand{\np}[1]{(#1)}
\newcommand{\bp}[1]{\big(#1\big)}

\newcommand{\horizon}{T}
\newcommand{\mos}{m}
\newcommand{\hum}{h}
\newcommand{\MOS}{M}
\newcommand{\HUM}{H}
\newcommand{\controlm}{u}
\newcommand{\dynamics}{\Phi}

\renewcommand{\SS}{\mathbb{S}}
\newcommand{\strategy}{\mathfrak{u}}

\newcommand{\robusto}{\mathbb{V}}

\newcommand{\RR}{\mathbb{R}}
\newcommand{\lb}[1]{\underline{#1}}  % {#1^{\flat}}
\newcommand{\ub}[1]{\overline{#1}}  % {#1^{\flat}}

\newcommand{\SSH}{\ub{\ub{\SS}}}
\newcommand{\SSM}{\ub{\SS}}
\newcommand{\SSL}{\mathring{\SS}}

\title{Robust Viability Analysis \\ of a 
Controlled Epidemiological Model}

\newcommand{\email}[1]{#1}
\newcommand{\keywords}[1]{\textbf{Keywords:} #1}

\author{Lilian Sofia Sepulveda Salcedo\footnote{%
Universidad Aut\'onoma de Occidente, 
Km. 3 v\'ia Cali-Jamund\'i, Cali, Colombia,
\email{lssepulveda@uao.edu.co}}
 \and 
Michel De Lara\footnote{%
Universit\'e Paris-Est, Cermics (ENPC), 
F-77455 Marne-la-Vall\'ee, %France,
\email{delara@cermics.enpc.fr}}   
}

\date{\today}

\begin{document}
\maketitle

\begin{abstract}
Managing infectious diseases is a world public health issue,
plagued by uncertainties.
In this paper, we analyze the problem of viable control of a dengue
outbreak under uncertainty.
For this purpose, we develop a controlled Ross-Macdonald model 
with mosquito vector control by fumigation,
and with uncertainties affecting the dynamics;
both controls and uncertainties are supposed to change only once a day, 
then remain stationary during the day.
The robust viability kernel is the set of all initial states such that 
there exists at least a strategy of insecticide spraying 
which guarantees that the number
of infected individuals remains below a threshold, 
for all times, and whatever the
sequences of uncertainties.
Having chosen three nested subsets of uncertainties --- 
a deterministic one (without uncertainty), a medium one and a large one ---
we can measure the incidence of the uncertainties
on the size of the kernel, in particular on its reduction with respect
to the deterministic case.
The numerical results show that the viability kernel without uncertainties 
is highly sensitive to the variability of parameters --- here the 
biting rate, the probability of infection to mosquitoes and humans,
and the proportion of female mosquitoes per person.
So the robust viability kernel is a possible tool to reveal the importance
of uncertainties regarding epidemics control. 
\end{abstract}


\pagebreak 
\tableofcontents

 \section{Introduction}

Managing infectious diseases is a world public health issue. 
The joint dynamics of infectious agents, vectors and hosts, 
and their spatial movements make the control or eradication problems difficult.
On top of that, uncertainties abound. 
Despite advances in epidemiological surveillance systems, 
and data on infected, recovered people, etc., there remain 
inaccuracies and errors. 
Factors such as ambient temperature, host age, social customs 
also contribute to uncertainty. 
In this paper, we focus on the impact of uncertainty on the viable control
of a Ross-Macdonald epidemiological model.

Our approach departs from the widespread approach in
mathematical epidemiology, where epidemic control aims at driving
the number of infected humans to zero, asymptotically.
Indeed, many studies on mathematical modeling of infectious diseases 
consist of analyzing the stability of the equilibria of a  
differential system (behavioral models such as SIR, SIS, SEIR
\cite{BrauerCastillo2006}). 
Those studies focus on asymptotic behavior and stability, 
generally leaving aside the transient behavior of the system, 
where the infection can reach high levels. 
In  many epidemiological models, a significant quantity is 
the ``basic reproductive number'' ${\cal R}_0$ 
which depends on parameters such as the transmission rate,
the mortality and birth rate, etc. 
Numerous works (see references in 
\cite{Hethcote:2000,Diekmann-Heesterbeek:2000})
exhibit conditions on ${\cal R}_0$ such that
the number of infected individuals tends towards zero. 
With this tool, different (time-stationary) 
management strategies of the propagation 
of the infection -- quarantine, vaccination, etc. -- are compared 
with respect to how they modify ${\cal R}_0$, that is, 
with respect to their capacity to drive 
the number of infected individuals towards zero, focusing on asymptotics.
However, during the transitory phase, the 
number of infected can peak at high values. 

By contrast, we focus both on the transitory and the asymptotic
regimes, where we aim at avoiding that the number of infected individuals
peaks at high values. As a consequence, our approach makes no reference the
concept of basic reproductive number~${\cal R}_0$. 

In \cite{DeLara-Sepulveda:2016}, we used viability theory to analyze 
the problem of maintaining the number of infected individuals
below a threshold, 
with limited fumigation capacity. The setting was deterministic, without
any uncertainty in the model parameters.
We said that a state is viable if there exists at least one 
admissible control trajectory 
--- time-dependent mosquito~mortality rates bounded by control~capacity --- 
such that, starting from this state, the resulting proportion of infected 
individuals remains below a given infection cap for all times.
We defined the so-called viability kernel as the set of viable states.
We obtained three different expressions of the viability kernel,
depending on the couple control~capacity-infection cap.

Several studies have applied the deterministic viable control method to 
managing natural resources, for example \cite{Bene-Doyen-Gabay:2001}, 
\cite{Bonneuil-Saint-Pierre:2005} and \cite{Bonneuil-Mullers:1997};
different examples can be found in \cite{DeLara-Doyen:2008}
in the discrete time case.
Yet, few studies have undertaken a robust approach to these issues
\cite{Bene-Doyen:2008,Regnier-DeLara:2015}. 

In this paper, we analyze what happens to the viability kernel when 
additional uncertain factors, that affect the dynamics of the disease, 
are considered. We make use of the so-called \emph{robust viability} approach
\cite{DeLara-Doyen:2008}, looking for vector fumigation policies able
to maintain the proportion of infected individuals below a 
given infection cap, for all times and for all scenarios of uncertainties.

Comparison of deterministic and robust viable states 
can contribute to shed light on the distance between the outcomes of these 
two extreme approaches: ignoring uncertainty versus hedging against any risk
(totally risk averse context \cite{Bene-Doyen:2008}).
In robust viability, constraints must be satisfied even under uncertainties 
related to an unlikely pessimistic scenario. 
By contrast, reducing uncertainties to zero amounts at dressing the problem as deterministic \cite{Aubin1991}.

The paper is organized as follows. 
In Section~\ref{The_Robust_Viability_Problem}, we introduce the 
controlled Ross-Macdonald model, discuss uncertainties and set the 
robust viability problem. Then, we introduce the robust viability kernel
and present the dynamic programming equation to compute it. 
In Section~\ref{sec:numeric calculations}, we provide a numerical application 
in the case of the dengue outbreak in 2013 in Cali, Colombia.  
We evaluate the impact of different sets of uncertainties on the 
robust viability kernel. We conclude in Section~\ref{sec:conclusions}.

\section{The robust viability problem}
\label{The_Robust_Viability_Problem}

In~\S\ref{sec:model_and_uncertainties}, we introduce the Ross-Macdonald 
model, we discuss uncertainties and we set the 
robust viability problem. Then, in~\S\ref{The_robust_viability_kernel}
we introduce the robust viability kernel.
In~\S\ref{Dynamic_programming_equation}, we present the dynamic programming 
equation to compute the robust viability kernel.

\subsection{Ross-Macdonald model with uncertainties}
\label{sec:model_and_uncertainties}

To address the robust viability problem, we work with a 
controlled Ross-Macdonald model with uncertainties. We suppose 
that both controls and uncertainties change values 
only at discrete time steps, then remain stationary between two consecutive 
time steps.

\subsubsection*{Controlled Ross-Macdonald model}
%\label{Ross-Macdonald_model_in_continuous_time}

Different types of Ross-Macdonald models have been published 
\cite{Smith_et_al2007}. We choose the one in \cite{AndersonMay1992},
where both total populations (humans, mosquitoes) are normalized to~1
and divided between susceptibles and infected. 
The basic assumptions of the model are the following.
\begin{itemize}
  \item[i)] The human population ($N_{\hum}$) and the mosquito population
($N_{\mos}$) are closed and remain stationary.
  \item[ii)] Humans and mosquitoes are homogeneous in terms of 
susceptibility, exposure and contact.
  \item[iii)] The incubation period is ignored, in humans as in mosquitoes.
  \item[iv)] Mortality induced by the disease is ignored, 
in humans as in mosquitoes.
  \item[v)] Once infected, mosquitoes never recover.
  \item[vi)] Only susceptibles get infected, in humans as in mosquitoes.
  \item[vii)] Gradual immunity in humans is ignored.
\end{itemize}

Time is continuous, denoted by $s \in \RR_+$.
Let $\mos(s) \in [0,1]$ denote the \emph{proportion of infected mosquitoes} 
at time~$s$,
and $\hum(s) \in [0,1]$ the \emph{proportion of infected humans} at time~$s$.
Therefore, $1-\mos(s)$ and $1-\hum(s)$ are the respective proportions 
of susceptibles.

The Ross-Macdonald model, used in \cite{DeLara-Sepulveda:2016},
is the following differential system
\begin{subequations}
\label{eq:model-Ross-Macdonald}
\begin{align}
  \displaystyle\frac{d\mos}{ds}&=
\alpha p_{\mos} \hum(1-\mos)-\delta \mos \eqfinv \\[5mm]
  \displaystyle\frac{d\hum}{ds}&=
\alpha p_{\hum}  \frac{N_{\mos}}{N_{\hum}}\, \mos (1-\hum)-\gamma \hum \eqfinv
\end{align}
\end{subequations}
where the parameters $\alpha$, $p_{\mos}$, $p_{\hum}$, 
$\xi=\frac{N_{\mos}}{N_{\hum}}$, 
$\delta$ and $\gamma$ are given in Table~\ref{tabla:parametros-Ross}.
The state space is the unit square~$[0,1]^2$.
\begin{table}[h]
	\centering
	\begin{tabular}{|l|l|l|}\hline
		{Parameter} & {Description} & {Unit} \\ \hline \hline
$ \alpha \geq 0 $ & biting rate per time unit & {time}$^{-1}$ \\ \hline
$  \xi=N_{\mos}/N_{\hum} \geq 0 $  & number of female mosquitoes per human & dimensionless
\\ \hline
$ 1 \geq p_{\hum} \geq 0 $ & probability of infection of a susceptible & \\
& human by infected mosquito biting & dimensionless \\ \hline
$ 1 \geq p_{\mos} \geq 0 $ & 
probability of infection of a susceptible  & \\
& mosquito when biting an infected human & dimensionless \\ \hline
$\gamma \geq 0 $  & recovery rate for humans & {time}$^{-1}$ \\ \hline
$\delta \geq 0 $  & (natural) mortality rate for mosquitoes & {time}$^{-1}$ \\ \hline
\end{tabular}
\caption{Parameters of the Ross-Macdonald model~\eqref{eq:model-Ross-Macdonald}}
\label{tabla:parametros-Ross}
\end{table}

For notational simplicity, we put
\begin{equation}
\label{eq:transmision-ross}
 A_{\mos}=\alpha p_{\mos} \eqsepv A_{\hum}=\alpha p_{\hum} \frac{N_{\mos}}{N_{\hum}}
=\alpha p_{\hum} \xi \eqfinp
\end{equation}
We turn the dynamical system~\eqref{eq:model-Ross-Macdonald} into a 
controlled system by replacing the natural mortality rate~$\delta$ 
for mosquitoes in~\eqref{eq:model-Ross-Macdonald} 
by a piecewise continuous function, called \emph{control trajectory}
\begin{equation}
  \controlm(\cdot): \RR_+ \to \RR \eqsepv s \mapsto \controlm(s) \eqfinp 
\label{eq:control_trajectory}
\end{equation}
Therefore, the \emph{controlled Ross-Macdonald model} is 
\begin{subequations}
  \begin{align}
  \displaystyle\frac{d\mos}{ds}=& 
A_{\mos}\hum(s) \bp{1-\mos(s)}-\controlm(s) \mos(s)  \eqsepv \\[5mm]
  \displaystyle\frac{d\hum}{ds}=& 
A_{\hum} \mos(s) \bp{1-\hum(s)}-\gamma \hum(s) \eqfinp   
  \end{align}
\label{eq:model-control-vector}
\end{subequations}
As is easily seen, the state space~$[0,1]^2$ is invariant by 
the dynamics~\eqref{eq:model-control-vector}.

\subsubsection*{Controlled Ross-Macdonald model sampled at discrete time steps}

Time steps, separated by a period of one day, are denoted by
\begin{equation}
  t = t_{0}, t_{0}+1, \dots, \horizon-1, \horizon \eqfinv
\label{eq:time}
\end{equation}
where $t_{0}\in\mathbb{N}$ is the \emph{initial time} (day)
and $\horizon\in\mathbb{N}$, $\horizon \geq t_{0}+1$ is the \emph{horizon}. 
Any interval $[t,t+1[$ represents one day.
Working with continuous time controls and uncertainties 
would make the mathematical framework more delicate, 
as well as the numerical calculation of the robust viability kernel. 
Thus our approach will be a mix, where the state follows a differential equation
in which controls and uncertainties remain stationary between two consecutive 
time steps.

From now on, we suppose that mosquito mortality rates induced by fumigation 
remain stationary during every time period (day), that is,
\begin{equation}
\forall t = t_{0}, t_{0}+1, \dots, \horizon-1 \eqsepv 
  \controlm(s)=\controlm(s') \eqsepv \forall \{s,s'\} \subset [t,t+1[
\eqfinp
\end{equation}
Let us denote by \( \dynamics(\MOS,\HUM, \controlm, A_{\MOS}, A_{\HUM}) \)
the solution, at time~$s=1$, of 
the differential system~\eqref{eq:model-control-vector}
with initial condition \( \bp{\mos(0),\hum(0)}=\np{\MOS,\HUM} \).
We obtain the following 
\emph{sampled and controlled Ross-Macdonald model}
\begin{equation}
\bp{\MOS(t+1),\HUM(t+1)} =
\dynamics\bp{\MOS(t),\HUM(t), \controlm(t), A_{\MOS}(t), A_{\HUM}(t)}
\label{eq:ross-incertidumbre}
%\label{eq:dinamica-sistema}
\end{equation}
where
\begin{itemize}
\item time~$t$ runs from~$t_{0}$ to ~$\horizon-1$, 
with time step one day, as in~\eqref{eq:time},
\item the state vector $\bp{\MOS(t), \HUM(t)} \in [0,1]^2$ 
represents the proportion of  mosquitoes and the proportion of humans 
that are infected at the beginning of day~$t$,
\item the control variable $\controlm(t)$ represents the mosquito mortality rate 
applied during all day~$[t,t+1[$, due to application of chemical control 
over the mosquitoes,
% \item the parameter $\gamma \geq 0$  is the recovery rate for humans 
% ({day}$^{-1}$), 
\item the vector $\big (A_{\MOS}(t), A_{\HUM}(t)\big ) \in \RR_+^2$ --- 
whose components are infection rates for mosquito and human, respectively
--- carries all the uncertainties; we discuss them below.
\end{itemize}

\subsubsection*{Uncertainties}

Uncertain factors in the dynamics of dengue transmission are 
the number of larval breeding sites, the rapid growth and urbanization 
of humans population, seasonal variability of mosquitoes population 
correlated with environmental factors such as rainfall, etc. 
The interactions between all these factors yield significant differences in 
parameters values --- such as mosquito bite rate, 
the proportion of mosquitoes females, 
probabilities of infection in humans and mosquitoes --- 
involved in the formulation of compartmental models of the process of 
transmission of an infectious diseases as dengue. 

In the Ross-Macdonald model~\eqref{eq:model-Ross-Macdonald},
we consider the following parameters as uncertain: 
bite rate ($\alpha$), 
probability of mosquito infection ($p_{\mos}$), 
probability of infection of a human ($p_{\HUM}$) and
number of female mosquitoes per person ($\frac{N_{\mos}}{N_{\hum}}$). 

As detailed in Appendix~\ref{Ranges_for_the_uncertain_aggregate_parameters},
we have decided to incorporate uncertainty in~\eqref{eq:ross-incertidumbre}
through the aggregate parameters 
$A_{\mos}=\alpha p_{\mos}$ and $A_{\hum}=\alpha p_{\hum} \frac{M}{H}$
in~\eqref{eq:transmision-ross}. Thus, the couples
\begin{equation}
\big (A_{\MOS}(t), A_{\HUM}(t)\big )\in \RR_+^2 \eqsepv
\forall t= t_0,\ldots, \horizon-1 
\end{equation}
represent uncertainty variables that affect the population of 
mosquitoes and humans population, respectively.

We define a \emph{scenario} of uncertainties as 
a sequence, of length~$\horizon-t_0$ of uncertainty couples:
\begin{equation}
\label{eq:escenario}
\begin{split}
\bp{A_{\MOS}(\cdot), A_{\HUM}(\cdot)}=  \\ 
\Big( \bp{A_{\MOS}(t_0),A_{\HUM}(t_0)}, \ldots, 
\bp{A_{\MOS}(\horizon-1),A_{\HUM}(\horizon-1)}\Big )
\in \np{\RR_+^2}^{\horizon-t_0} \eqfinp
\end{split}
\end{equation}

\subsection{The robust viability kernel}
\label{The_robust_viability_kernel}

Viable management in the robust sense evokes a pessimistic context. 
Here, it will mean that we want to satisfy the constraint of keeping 
the number of infected humans below a threshold imposed, 
\emph{whatever the uncertainties}. We now give precise mathematical statement
of this problem.

\subsubsection*{Control constraints}

Let \( (\lb{\controlm}, \ub{\controlm}) \) be a couple such that 
\begin{equation}
0 \leq \lb{\controlm} \leq \ub{\controlm} \leq 1 \eqfinp
\label{eq:upper_and_lower_control_bounds}
\end{equation}
We impose that the control variable~$\controlm(t)$ --- 
the mosquito mortality rate applied during all day~$[t,t+1[$ 
in~\eqref{eq:ross-incertidumbre} ---
satisfies the constraints
\begin{equation}
 \lb{\controlm} \leq \controlm(t) \leq \ub{\controlm} \eqsepv
\forall t= t_0,\ldots, \horizon-1 \eqfinp
\label{eq:restriccion-control-incertidumbre}
\end{equation}

\subsubsection*{Uncertainty constraints}

The two terms $A_{\MOS}(t)$ and $A_{\HUM}(t)$ in~\eqref{eq:ross-incertidumbre} 
encapsulate all the uncertainties affecting each population, respectively. 
We suppose that they can take any value in a known 
\emph{set~$\SS$ of uncertainties}: 
\begin{equation}
\bp{A_{\MOS}(t), A_{\HUM}(t)}\in \SS \subset \RR_+^2 \eqsepv
\forall t= t_0,\ldots, \horizon-1 \eqfinp
\label{eq:set_of_uncertainties}
\end{equation}
We could have allowed the set~$\SS$ to differ from one time to another
(hence being a sequence of $\SS(t)$ for $t= t_0,\ldots, \horizon-1$).
However, for simplicity reasons, we will only consider the stationary case. 

\subsubsection*{State constraints}

Let~$\ub{\HUM}$ be a real number such that 
\begin{equation}
  0< \ub{\HUM} < 1 \eqfinv
\label{eq:cap}
\end{equation}
which represents the \emph{maximum tolerated proportion of infected humans},
or \emph{infection cap}. 
Thinking about public health policies set by governmental entities, 
we impose the following constraint:
the proportion~$\HUM(t)$ of infected humans 
must always remain below the {infection cap}~$\ub{\HUM}$.
Therefore, we impose the state constraint
\begin{equation}
\label{eq:restriccion-estado-incertidumbre}
\HUM(t) \leq \ub{\HUM} \eqsepv \forall t=t_0,\ldots, \horizon-1,T \eqfinp
\end{equation}
For example, the Municipal Secretariat of Public Health of Cali, Colombia,
establishes a so-called ``endemic canal'' (\emph{canal end\'emico}) 
with an upper bound on infected individuals.
 
\subsubsection*{Strategies}

To define the robust viability kernel, we need the notion of strategy.
A  \emph{control strategy}~$\strategy$ is a sequence of mappings 
from states~$(\MOS, \HUM)$ towards controls~$\controlm$ as follows:
\begin{equation}
\label{eq:estrategia}
\strategy=\np{\strategy_t}_{t=t_0,\ldots, \horizon-1} \text{ with } 
\strategy_t:[0,1]^2 \to \RR \eqfinp 
\end{equation}
A control strategy~$\strategy$ is said to be \emph{admissible} if 
\begin{equation}
\label{eq:estrategia-admisible}
\strategy_t:[0,1]^2 \to [\lb{\controlm}, \ub{\controlm}] \eqfinp 
\end{equation}
\begin{subequations}
For any scenario \( \bp{A_{\MOS}(\cdot), A_{\HUM}(\cdot)} \) 
as in~\eqref{eq:escenario}, and any initial state
\begin{equation}
\bp{\MOS(0), \HUM(0)}= (\MOS_0, \HUM_0) \in [0,1]^2 \eqfinv
\label{eq:initial_state}
\end{equation}
a control strategy~$\strategy$ as in~\eqref{eq:estrategia} produces
--- through the dynamics~\eqref{eq:ross-incertidumbre} --- 
a state trajectory~$\bp{\MOS(\cdot), \HUM(\cdot)}$ by the closed-loop dynamics
\begin{equation}
\bp{\MOS(t+1),\HUM(t+1)} =
\dynamics\bp{\MOS(t),\HUM(t), \strategy_t\bp{\MOS(t), \HUM(t)},
A_{\MOS}(t), A_{\HUM}(t)} \eqfinv
\label{eq:dinamica-lazo-cerrado}
\end{equation}
for $t=t_0, \ldots, \horizon-1$, 
and a control trajectory~$\controlm(\cdot)$ by
\begin{equation}
\label{eq:trayectoria-control}
\controlm(t)=\strategy_t\bp{\MOS(t), \HUM(t)}\eqsepv 
t=t_0, \ldots, \horizon-1 \eqfinp
\end{equation}
\label{eq:closed-loop}
\end{subequations}

\subsubsection*{Robust viability  kernel}

Now, we are ready to lay out the definition of the robust viability kernel.

\begin{definition}
The \emph{robust viability kernel} (at initial time~$t_0$) is
\begin{equation}
\robusto_{\SS}(t_0)=
\left \{
\begin{array}{c}
(\MOS_0, \HUM_0) \\ \in [0,1]^2 
\end{array} 
\left |
\begin{array}{c}
 \textit{there is at least one} \\
 \textit{admissible strategy } \strategy 
 \textit{ as in }~\eqref{eq:estrategia-admisible}\\
 \textit{ such that, for every scenario } \\
\bp{A_{\MOS}(\cdot), A_{\HUM}(\cdot)} \in \SS^{\horizon-t_0}, \\
  \textit{the state trajectory } \\
\bp{\MOS(\cdot), \HUM(\cdot)}
\textit{ given by~\eqref{eq:closed-loop} }\\
  \textit{satisfies the state constraint }~\eqref{eq:restriccion-estado-incertidumbre} \\     
\end{array} \right.
\right \}\eqfinp
\label{eq:nucleo-viabilidad-robusto}
\end{equation}
\label{de:nucleo-viabilidad-robusto}
\end{definition}
The states belonging to the robust viability kernel are called 
\emph{viable robust states}. 

We use the notation~$\robusto_{\SS}(t_0)$ 
to mark the dependency of the robust viability kernel
of~\eqref{eq:nucleo-viabilidad-robusto} with respect to the set~$\SS$
of uncertainties in~\eqref{eq:set_of_uncertainties}. 
Indeed, in the sequel, we will study how $\robusto_{\SS}(t_0)$ varies with~$\SS$.

\subsection{Dynamic programming equation}
\label{Dynamic_programming_equation}

In~\eqref{eq:nucleo-viabilidad-robusto}, we observe that 
the set of scenarios --- with respect to which 
the robust viability kernel is defined --- 
is the rectangle \( \SS^{\horizon-t_0} \).
This is a strong assumption that corresponds to stagewise independence 
of uncertainties. Indeed, when the ``head'' uncertainty trajectory
\( \Big( \bp{A_{\MOS}(t_0),A_{\HUM}(t_0)}, \ldots,\)
\( \bp{A_{\MOS}(t),A_{\HUM}(t}\Big ) \) is known,
the domain where the ``tail'' uncertainty trajectory
\( \Big( \bp{A_{\MOS}(t+1),A_{\HUM}(t+1)}, \ldots,
\bp{A_{\MOS}(\horizon-1),A_{\HUM}(\horizon-1)}\Big ) \)
takes its value does not depend on the head uncertainty trajectory,
as it is the rectangle \( \SS^{\horizon-t+1} \).

This rectangularity property of the set of scenarios makes it possible
to compute the robust viability kernel~$\robusto_{\SS}(t_0)$ 
of~\eqref{eq:nucleo-viabilidad-robusto} by dynamic programming.
The proof of the following result is an easy extension of a proof
to be found in~\cite{DeLara-Doyen:2008}.

Let us introduce the so-called \emph{constraints set} 
\begin{equation}
\mathbb{A}=\{(\MOS,\HUM)|0\leq\MOS\leq 1 \eqsepv 0\leq\HUM\leq\ub{\HUM}\}
= [0,1] \times [0,\ub{\HUM}] 
\end{equation}
and let $1_{\mathbb{A}}(\MOS,\HUM)$ denote its \emph{indicator function} 
\begin{equation}
1_{\mathbb{A}}(\MOS,\HUM)=
\Bigg \{
\begin{array}{cc}
1 & \text{if } (\MOS,\HUM)\in\mathbb{A} \eqfinv \\
0 & \text{if } (\MOS,\HUM)\notin\mathbb{A} \eqfinp
\end{array}
\end{equation}

\begin{proposition}
The robust viability kernel~$\robusto_{\SS}(t_0)$ 
of~\eqref{eq:nucleo-viabilidad-robusto} is given by 
\begin{equation}
  \robusto_{\SS}(t_0) = \{ (\MOS,\HUM) \in [0,1]^2 \mid 
\mathbf{V}_{t_0}(\MOS,\HUM) =1 \} \eqfinv
\end{equation}
where the function $\mathbf{V}_{t_0}$ is solution of the 
following backward induction
--- called \emph{dynamic programming equation} --- 
that connects the so-called \emph{value functions}
\begin{subequations}
\label{eq:programacion-dinamica}
\begin{align}
\mathbf{V}_T(\MOS,\HUM) &= 1_{\mathbb{A}}(\MOS,\HUM)
\eqsepv \forall (\MOS,\HUM) \in [0,1]^2 \eqfinv
\\[5mm]
\mathbf{V}_t(\MOS,\HUM) &= 1_{\mathbb{A}}(\MOS,\HUM) 
\displaystyle\sup_{\controlm\in[\lb{\controlm},\ub{\controlm}]}
\displaystyle\inf_{(A_{\MOS}, A_{\HUM})\in\SS}
\mathbf{V}_{t+1}\bp{\dynamics(\MOS,\HUM, \controlm, A_{\MOS}, A_{\HUM})}
\eqsepv \\
& \qquad \forall (\MOS,\HUM) \in [0,1]^2 \eqfinv
\end{align}
where $t$ runs down from $\horizon-1$ to $t_0$.
\end{subequations} 
\end{proposition}

We are now equipped with the concept of robust viability kernel
and we dispose of a method to compute it.

\section{Numerical results and viability analysis}
\label{sec:numeric calculations}

Now, after having set the theory in 
Section~\ref{The_Robust_Viability_Problem}, 
we will provide a numerical application 
in the case of the dengue outbreak in 2013 in Cali, Colombia.  
We introduce three nested sets of uncertainties
in~\S\ref{Three_nested_sets_of_uncertainties},
we study their impact on the robust viability kernels computed
in~\S\ref{Numerical_computation_of_robust_viability_kernels},
and, finally, we discuss in~\S\ref{Robust_viability_analysis}
the results thus obtained.

\subsection{Three nested sets of uncertainties}
\label{Three_nested_sets_of_uncertainties}

With the purpose of evaluating the sensitivity of the size and shape 
of the robust viability kernel~$\robusto_{\SS}(t_0)$ 
of~\eqref{eq:nucleo-viabilidad-robusto} with respect to the set~$\SS$
of uncertainties in~\eqref{eq:set_of_uncertainties}, 
we consider the following three possible cases --- 
$\SSL$, $\SSM$ and $\SSH$ --- for the set~$\SS$:
\begin{subequations}
\begin{description}
\item[L)] Low case (deterministic)
\begin{equation}
\SSL=\Big \{\mathring{A}_{\MOS}\Big \}\times \Big \{\mathring{A}_{\HUM}\Big \}
\eqfinv
\label{eq:determinista}
\end{equation}
\item[M)] Middle case
\begin{equation}
\SSM=\Big [\lb{A}_{\MOS}, \ub{A}_{\MOS}\Big ]\times\Big [\lb{A}_{\HUM}, \ub{A}_{\HUM}\Big ] \eqfinv
\label{eq:caso-media}
\end{equation}
\item[H)] High case
\begin{equation}
\SSH=\Big [\lb{\lb{A}}_{\MOS}, \ub{\ub{A}}_{\MOS}\Big ]\times\Big [\lb{\lb{A}}_{\HUM}, \ub{\ub{A}}_{\HUM}\Big ] \eqfinp
\label{eq:caso-alta}
\end{equation}
\end{description}
\end{subequations}
We make the additional assumption that
\begin{equation}
\SSL\subset\SSM\subset\SSH\subset \RR_+^2\eqfinv
\end{equation} 
that is,
\begin{subequations}
\begin{align}
0 \leq \lb{\lb{A}}_{\MOS} \leq \lb{A}_{\MOS} \leq \mathring{A}_{\MOS} 
\leq \ub{A}_{\MOS} \leq \ub{\ub{A}}_{\MOS} \eqfinv \\
0 \leq \lb{\lb{A}}_{\HUM} \leq \lb{A}_{\HUM} \leq \mathring{A}_{\HUM} 
\leq \ub{A}_{\HUM} \leq \ub{\ub{A}}_{\HUM} \eqfinp
\end{align}
\end{subequations}
By Definition~\ref{de:nucleo-viabilidad-robusto}, % --- 
% and with the notation $\robusto_{\SS}(t_0)$ to mark the dependency 
% of the robust viability kernel~$\robusto(t_0)$ 
% of~\eqref{eq:nucleo-viabilidad-robusto} with respect to the set~$\SS$
% of uncertainties in~\eqref{eq:set_of_uncertainties} --- 
 we easily obtain that
\begin{equation}
\robusto_{\SSH}(t_0)\subset\robusto_{\SSM}(t_0)\subset\robusto_{\SSL}(t_0)
\eqfinv
\end{equation}
as less and less initial states can comply with 
the constraints in~\eqref{eq:nucleo-viabilidad-robusto}
as the set~$\SS$ of uncertainties increases.

\subsection{Numerical computation of robust viability kernels}
\label{Numerical_computation_of_robust_viability_kernels}

We discuss the choice of the following numerical values 
in Appendix~\ref{Appendix_Epidemic_model}.
For the numerical calculation of the robust viability kernel
in~\eqref{eq:nucleo-viabilidad-robusto}, 
we use the dynamic programming equation~\eqref{eq:programacion-dinamica}
with over 60~days, that is, 
\begin{equation}
t_{0}=0 \eqsepv \horizon= 60~\text{days.} 
\end{equation}
Following the recommendation of 
To remain below the ``endemic canal'' (\emph{canal end\'emico}) 
established by the Municipal Secretariat of Public Health 
(Cali, Colombia), we take for the infection cap the value
\begin{equation}
\ub{\HUM} = 0.00001 = 0.001\% \eqfinp   
\end{equation}
With this value, we express that we do not want that more
than 0.001\% of the total population of humans be infected. 
We will also increase this value to \( \ub{\HUM} = 0.0001 = 0.01\% \),
to test how the robust viability kernel is impacted.

We discretize the state space~$[0, 1]^2$, the control space
(see~\eqref{eq:natural_mortality_rate_estimated}-%
\eqref{eq:maximal_mosquito_mortality_rate_estimated}
for the choice of numerical values)
\begin{equation}
  [\lb{\controlm},\ub{\controlm}] = [0.0333,~0.05]~\textrm{day}^{-1}
\end{equation}
and the sets~$\SSM$ 
and $\SSH$ of uncertainties in~\eqref{eq:caso-media}-\eqref{eq:caso-alta}. 
We take a partition of $70$~elements for
\begin{itemize}
\item 
the interval $[0, 1]$, where the proportion~$\MOS$
of infected mosquitoes takes its values,
\item 
the interval $[0, \ub{\HUM}]$, where the proportion~$\HUM$
of infected humans takes its values,
\item 
the interval $[\lb{\controlm},\ub{\controlm}]$, 
where the control variable~$\controlm$  takes its values,
\item 
the intervals $[\lb{A}_{\MOS}, \ub{A}_{\MOS}]$ and $[\lb{A}_{\HUM}, \ub{A}_{\HUM}]$
(or the intervals $[\lb{\lb{A}}_{\MOS}, \ub{\ub{A}}_{\MOS}]$ 
and $[\lb{\lb{A}}_{\HUM}, \ub{\ub{A}}_{\HUM}]$)
where the uncertainties ${A}_{\MOS}(t)$ and ${A}_{\HUM}(t)$ 
take their values, respectively, as 
 in~\eqref{eq:caso-media}-\eqref{eq:caso-alta}.
\end{itemize}
To compute the robust viability kernel~$\robusto_{\SS}(t_0)$
by~\eqref{eq:programacion-dinamica}, 
we implement the following dynamic programming 
Algorithm~\ref{algoritmo-nucleo-robusto}. 
\bigskip

\begin{algorithm}[H]
\vspace{2mm}

 Initialization 
$ \mathbf{V}_T(\MOS,\HUM) = 1_{\mathbb{A}}(\MOS,\HUM) $\;\vspace{2mm}

\For{ $t=\horizon-1, \dots, t_0$ }{\vspace{2mm}

\ForAll { $(\MOS,\HUM)$ }{\vspace{2mm}

%\ForAll{ $\controlm  \in [\lb{\controlm},\ub{\controlm}]$ }{\vspace{2mm}
\ForAll{ $\controlm $ }{\vspace{2mm}

% \ForAll{ $(A_{\MOS}, A_{\HUM})\in\SS$ }{\vspace{2mm}
\ForAll{ $(A_{\MOS}, A_{\HUM})$ }{\vspace{2mm}

$ \mathbf{V}_{t+1}\bp{\dynamics(\MOS,\HUM,\controlm,A_{\MOS}, A_{\HUM})}$ }
\vspace{2mm}
% $ \displaystyle\min_{(A_{\MOS}, A_{\HUM})\in\SS} 
$ \displaystyle\min_{(A_{\MOS}, A_{\HUM})} 
\mathbf{V}_{t+1}\bp{\dynamics(\MOS,\HUM,\controlm,A_{\MOS}, A_{\HUM})} $ 
}\vspace{2mm}

% $ \displaystyle\max_{\controlm  \in [\lb{\controlm},\ub{\controlm}]} 
% \min_{(A_{\MOS}, A_{\HUM})\in\SS} 
$ \displaystyle\max_{\controlm } \min_{(A_{\MOS}, A_{\HUM})} 
\mathbf{V}_{t+1}\bp{\dynamics(\MOS,\HUM,\controlm,A_{\MOS}, A_{\HUM})} $ }
\vspace{2mm}
$ \mathbf{V}_{t}(\MOS,\HUM)=1_{\mathbb{A}}(\MOS,\HUM) 
\mathbf{V}_{t+1}\bp{\dynamics(\MOS,\HUM,\controlm,A_{\MOS}, A_{\HUM})}$
}
\vspace{2mm}
\caption{\label{algoritmo-nucleo-robusto} Dynamic programming algorithm
to compute the robust viability kernel by~\eqref{eq:programacion-dinamica}}
\end{algorithm}
\bigskip

In general, the image $\dynamics(\MOS,\HUM,\controlm,A_{\MOS}, A_{\HUM})$ 
of the dynamics does not fall exactly on one element of the $70\times 70$ grid 
--- corresponding to the cells of the state space 
$[0, 1]\times[0, \ub{\HUM}]$ --- over which the 
numerical value function~$\mathbf{V}_{t+1}$ is defined. 
We have taken a conservative stand: we put
\( \mathbf{V}_{t+1}\bp{\dynamics(\MOS,\HUM,\controlm,A_{\MOS}, A_{\HUM})} 
= 1 \) if and only if \( \mathbf{V}_{t+1}=1 \) for 
\emph{all} points in the grid that surround
the image $\dynamics(\MOS,\HUM,\controlm,A_{\MOS}, A_{\HUM})$. 

The (discrete) robust viability kernel is defined as the set of points 
of the $70\times 70$ grid where the indicator function 
$\mathbf{V}_{t_0}\bigl(\MOS,\HUM\bigr)$ is equal to~1.

\subsection{Robust viability analysis}
\label{Robust_viability_analysis}

Now, as we are able to compute robust viability kernels 
as just seen in~\S\ref{Numerical_computation_of_robust_viability_kernels},
we will compare them under the three nested sets of uncertainties
introduced in~\S\ref{Three_nested_sets_of_uncertainties}.

\subsubsection*{Comparison between deterministic viability kernels 
in discrete time and in continuous time}

We consider the low case (deterministic), for which we take  
\begin{equation}
\SSL= \{0.07660\}\times\{0.0722\}~\text{day}^{-1} \times \text{day}^{-1} 
\eqfinp
\end{equation}
These numbers were obtained in~\eqref{eq:transmision-ross_estimated}
from the parameters adjusted to the 2013 dengue outbreak in Cali, Colombia,
as described in~\cite{DeLara-Sepulveda:2016}. Details can be found
in Appendix~\ref{Appendix_Epidemic_model}.

\begin{figure}
\begin{center}
\includegraphics[width=13cm,height=7cm]{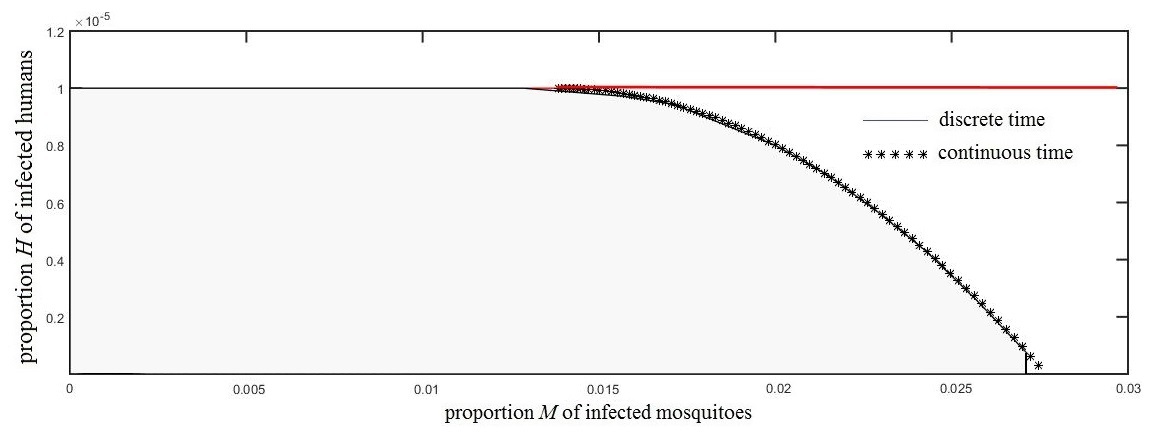}
\end{center}
\caption{Viability kernels for the deterministic case:
obtained by the dynamic programming Algorithm~\ref{algoritmo-nucleo-robusto}
in discrete time ($\SSL=\{0.07660\}\times\{0.0722\}$, dotted line) 
and given by a mathematical formula in~\cite{DeLara-Sepulveda:2016}
in continuous time (broken starred line),
for the infection cap \( \ub{\HUM} = 0.00001 = 0.001\% \)}
\label{fig:nucleo-determinista5}
\end{figure}

\begin{figure}
\begin{center}
\includegraphics[width=13cm,height=7cm]{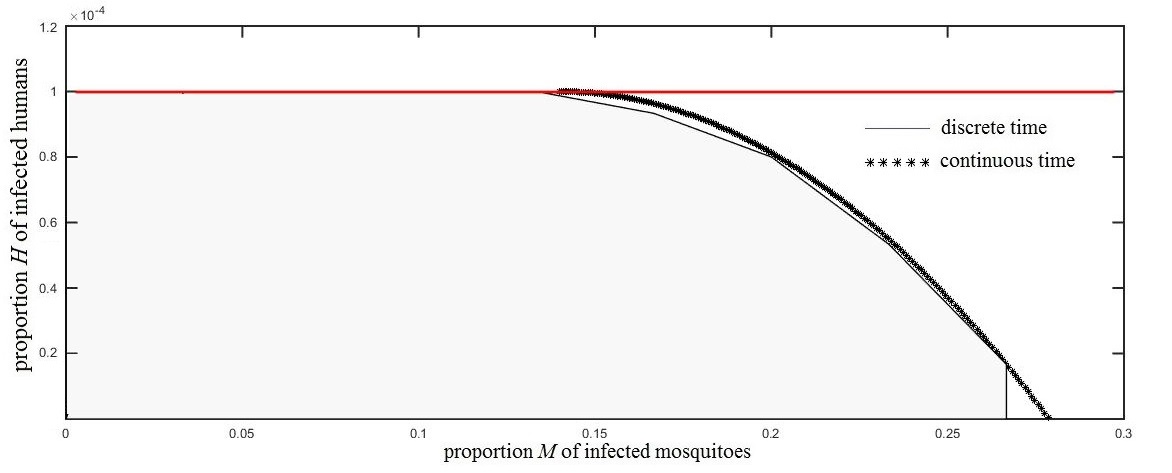}
\end{center}
\caption{Viability kernels for the deterministic case:
obtained by the dynamic programming Algorithm~\ref{algoritmo-nucleo-robusto}
in discrete time ($\SSL=\{0.07660\}\times\{0.0722\}$, dotted line) 
and given by a mathematical formula in~\cite{DeLara-Sepulveda:2016}
in continuous time (broken starred line),
for the infection cap \( \ub{\HUM} = 0.0001 = 0.01\% \)}
\label{fig:nucleo-determinista4}
\end{figure}

%In Figure~\ref{fig:nucleo-determinista5}, 
In Figures~\ref{fig:nucleo-determinista5} and~\ref{fig:nucleo-determinista4}, 
we present both the deterministic viability kernel in discrete time --- 
that is, obtained when $\SS=\SSL$ in~\eqref{eq:determinista} 
and computed with the dynamic programming 
Algorithm~\ref{algoritmo-nucleo-robusto}---
and the deterministic viability kernel in continuous time --- 
that is, the one given by a mathematical formula 
in~\cite{DeLara-Sepulveda:2016}. 
The solid continuous line corresponds to the discrete time case 
and the broken starred line corresponds to the continuous time case. 

We observe that the viability kernel 
for the deterministic case in discrete time is almost identical
to the deterministic viability kernel in continuous time.
Therefore, neither the sampling of controls and uncertainties --- as described 
in~\S\ref{sec:model_and_uncertainties} ---
nor the discretization process --- described 
in~\S\ref{Numerical_computation_of_robust_viability_kernels} ---
lead to a degradation of the theoretical viability kernel.

\subsubsection*{Comparison between robust viability kernels 
as the set~$\SS$ of uncertainties increases}

\begin{figure}
\centering
\includegraphics[width=13cm,height=7cm]{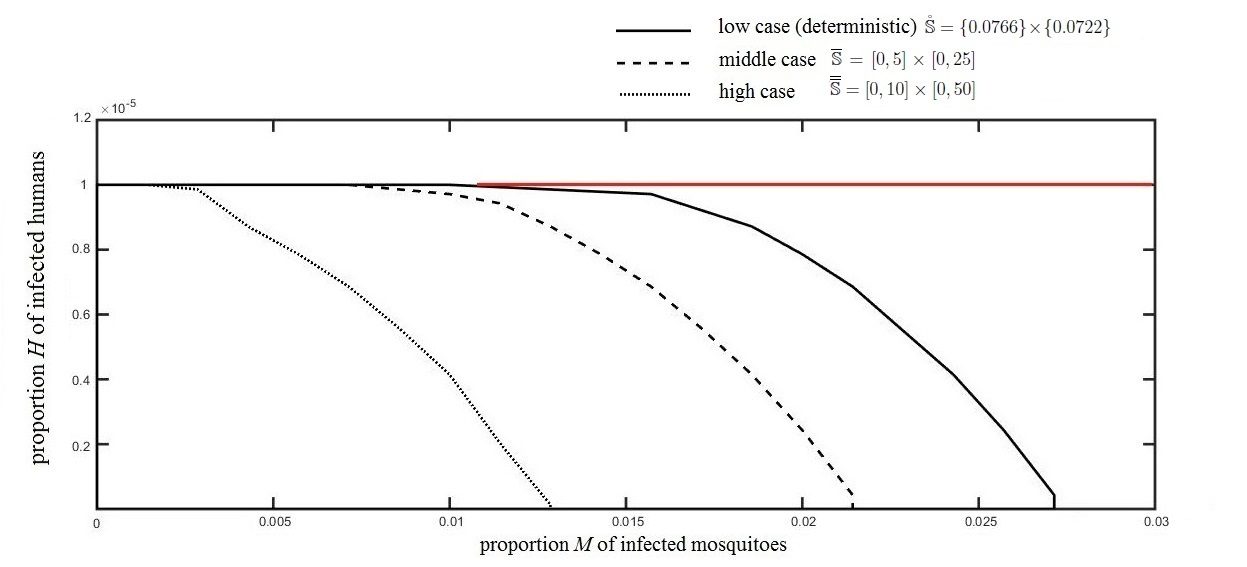}
\caption{Deterministic viability kernel 
($\SSL=\{0.0766\}\times\{0.0722\}$, right hand side continuous line) 
and robust viability kernel for the middle case 
($\SSM=[0,5]\times[0, 25]$, middle broken line) 
and the high case 
($\SSH=[0,10]\times[0, 50]$, left hand side dotted line),
for the infection cap \( \ub{\HUM} = 0.00001 = 0.001\% \)}
\label{fig:nucleo-robusto5}
\end{figure}

\begin{figure}
\centering
\includegraphics[width=13cm,height=7cm]{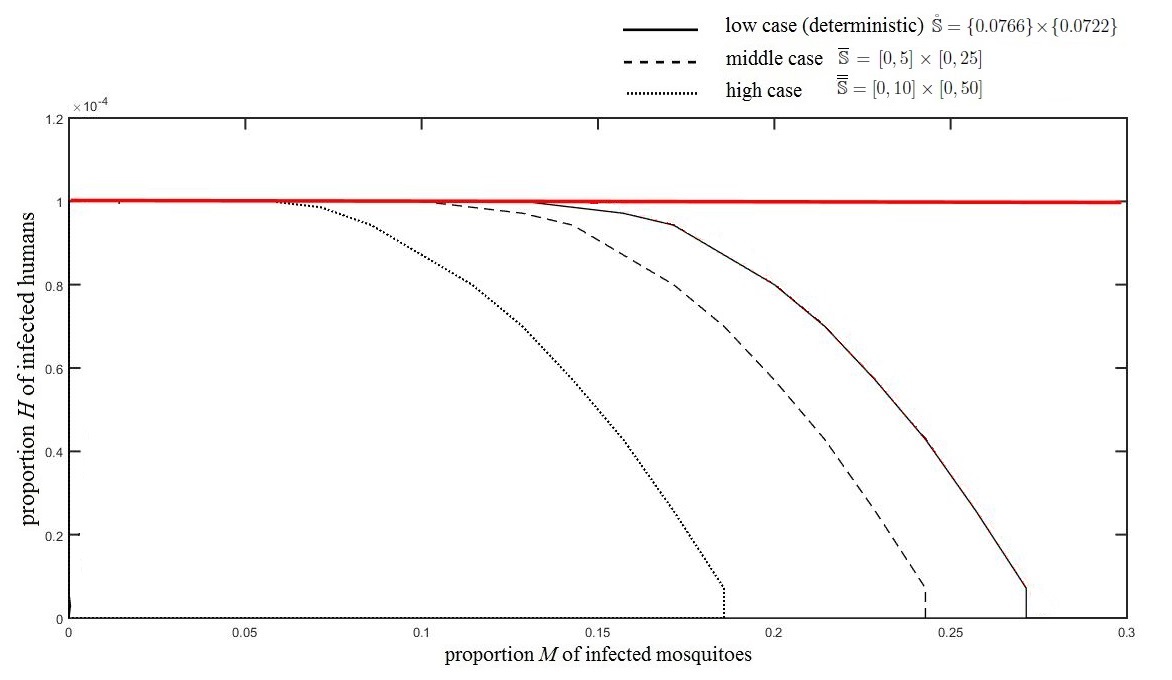}
\caption{Deterministic viability kernel 
($\SSL=\{0.0766\}\times\{0.0722\}$, right hand side continuous line) 
and robust viability kernel for the middle case 
($\SSM=[0,5]\times[0, 25]$, middle broken line) 
and the high case 
($\SSH=[0,10]\times[0, 50]$, left hand side dotted line),
for the infection cap \( \ub{\HUM} = 0.0001 = 0.01\% \)}
\label{fig:nucleo-robusto4}
\end{figure}

%In Figure~\ref{fig:nucleo-robusto}, 
In Figures~\ref{fig:nucleo-robusto5} and~\ref{fig:nucleo-robusto4}, 
we display the three robust viability kernels 
corresponding to the three cases presented 
in~\S\ref{Three_nested_sets_of_uncertainties}:
low case (deterministic) with~$\SS=\SSL$, middle case with~$\SS=\SSM$ and 
high case with~$\SS=\SSH$. 

For the middle case, we take  
\begin{equation}
\label{eq:incertidumbre-media}
\SSM=[0,5]\times[0, 25]~\text{day}^{-1} \times \text{day}^{-1} \eqfinp
\end{equation}
These numbers were obtained from the ranges for the parameters 
adjusted to the 2013 dengue outbreak in Cali, Colombia,
as described in~\cite{DeLara-Sepulveda:2016}. Details can be found
in Appendix~\ref{Ranges_for_the_uncertain_aggregate_parameters}.

For the high case, we doubled the right ends of each 
interval in~\eqref{eq:incertidumbre-media}, giving:
\begin{equation}
\label{eq:incertidumbre-alta}
\SSH=[0,10]\times[0, 50]~\text{day}^{-1} \times \text{day}^{-1} \eqfinp
\end{equation}


\subsubsection*{Discussion}

In Figures~\ref{fig:nucleo-robusto5} and~\ref{fig:nucleo-robusto4}
%In Figure~\ref{fig:nucleo-robusto} 
--- both in the middle case (middle dashed line)
and more markedly in the high case (left hand side dotted line) --- 
we observe the following:
a large part of the initial states~$(\MOS_0, \HUM_0)$ that are 
identified as viable in the deterministic case 
(right hand side continuous line)
are no longer viable when taking into account uncertainties.

In addition, as we expand the set~$\SS$ where uncertainties take their values, 
the gap with the deterministic case is larger and larger. 
In the high case, when we double the extremites of~$\SSM$, 
we observe in Figures~\ref{fig:nucleo-robusto5} and~\ref{fig:nucleo-robusto4}
% Figure~\ref{fig:nucleo-robusto} 
that all the initial states that are below the 
right hand side continuous line 
(boundary of the deterministic viability kernel) 
and above the lower left hand side dotted line
no longer belong to the viability kernel.
Thus, enlarging the set of uncertainties can have a strong impact 
on the viability kernel. When we have much variability in the values of 
the uncertain variables~$(A_{\mos}, A_{\HUM})$, 
the robust viability kernel is greatly reduced.

Therefore, we cannot guarantee compliance with the constraint 
imposed on the proportion of infected humans for all uncertainties, 
if we take initial conditions belonging to the deterministic kernel.
Picking initial conditions within the deterministic kernel
would yield optimistic conclusions with respect to infection control,
when uncertainties affect the epidemics dynamics.

These observations depend on the temporal horizon and on the
structure of the set of scenarios. Recall that, 
in~\S\ref{Dynamic_programming_equation},
we stressed the importance of having a rectangular set of scenarios 
to obtain a dynamic programming equation, making it possible
to recursively compute robust viability kernels. 
This rectangularity property of the set of scenarios 
corresponds to stagewise independence of uncertainties:
whatever the value of the uncertainty at time~$t$, 
the next uncertainties can take \emph{any} values 
(within the rectangle).
This rectangularity property 
contributes to having small robust viability kernels. 
Indeed, scenarios can display arbitrary evolutions, 
switching from one extreme to another between time periods (days). 
High and low values for the uncertainties alternate,
submitting the epidemics dynamics to a strong stress,
and thus narrowing the possibility of satisfying 
the state constraints for all times.
Such contrasted scenarios deserve the label of worst-case scenarios
(they represent a small fraction of the extreme points of the rectangular set of scenarios).
This is why amplifying the distance 
between extreme uncertainties shrinks the robust viability kernel.

\section{Conclusion}
\label{sec:conclusions}

In~\cite{DeLara-Sepulveda:2016}, the two authors obtained a neat 
expression of the deterministic viability kernel for a 
controlled Ross-Macdonald model.
However, uncertainties abound and we wanted to assess their importance 
regarding epidemics control. 
The numerical results show that the viability kernel without uncertainties 
is highly sensitive to the variability of parameters --- here the 
biting rate, the probability of infection to mosquitoes and humans,
and the proportion of female mosquitoes per person.
So, a robust viability analysis can be a tool to reveal the importance
of uncertainties regarding epidemics control. 

\appendix

\section{Appendix. Fitting an epidemiological model for dengue}
\label{Appendix_Epidemic_model}
\newcommand{\zz}{z}
\newcommand{\thetab}{\theta}

Here, we present how we identify parameters for the 
Ross-Macdonald model~\eqref{eq:model-Ross-Macdonald},
and then obtain ranges for the uncertain aggregate parameters. 

\subsection{Parameters and daily data deduced from health reports}

We introduce % the state vector 
% \begin{equation}
%   \zz= ( \mos, \hum) \in [0,1]^2 \eqfinv
% \end{equation}
% and 
the vector of parameters
\begin{equation}
\thetab=\big( \alpha, p_{\hum}, p_{\mos}, \xi, \delta \big) \in \Theta 
 \subset \RR^5_+  \eqfinv
\label{eq:parameter}
\end{equation}
consisting of the five parameters previously defined in 
Table~\ref{tabla:parametros-Ross}. The parameter set 
$\Theta \subset \RR^5_+$ is given by the Cartesian product of the 
five intervals in the third column of Table~\ref{tab-A1}. 

With these notations, the Ross-Macdonald model~\eqref{eq:model-Ross-Macdonald}
now writes
\begin{equation}
\label{RMmodel-par}
\begin{array}{rl}
  \displaystyle\frac{d\mos(s;\thetab)}{ds}&=
\alpha p_{\mos} \hum(s;\thetab) \bp{1-\mos(s;\thetab)}-\delta \mos
\eqfinv \\[5mm]
  \displaystyle\frac{d\hum(s;\thetab)}{ds}&=
\alpha p_{\hum} \xi \mos(s;\thetab) \bp{1-\hum(s;\thetab)}-\gamma \hum 
\eqfinv \\[5mm] 
\bp{\mos(s_0;\thetab),\hum(s_0;\thetab)} &= (\mos_0, \hum_0) \eqfinp
\end{array} % \right) \eqsepv \gamma=0.1 \eqfinp 
\end{equation}

Notice that the rate~$\gamma$ of human recovery 
does not appear in the parameter vector~$\thetab$ in~\eqref{eq:parameter}. 
Indeed, in the data provided by the Municipal Secretariat of Public Health 
(Cali, Colombia), we only have new cases of dengue registered per day; 
there is no information regarding how many inviduals recover daily. 
We choose an infectiousness period of~10 days, that is, 
a rate of human recovery fixed at 
\begin{equation}
  \gamma=0.1~\text{day}^{-1} \eqfinp
\label{eq:gamma_estimated}
\end{equation}
Under this assumption, the \emph{daily incidence data} 
(i.e., numbers of newly registered cases reported on daily basis) 
provided by the Municipal Secretariat of Public Health
can be converted into the \emph{daily prevalence data}
(i.e., numbers of infected inviduals on a given day, be they new or not).
With this, we deduce values of daily proportion of infected inviduals
in the form of the set
\begin{equation}
\label{dataset}
\mathbb{O} = \Big\{ \big( r, \hat{\hum}_r \big) \eqsepv 
r =r_0,r_0 +1, \ldots, R \Big\} \eqfinv 
\end{equation}
where $r$ refers to day~$r$ within the observation period 
of~$R+1-r_0$ days, and where $\hat{\hum}_r$ stands for the 
fraction of infected inviduals at day~$r$. Naturally, the first couple 
in the set~$\mathbb{O}$ in \eqref{dataset} defines the initial condition 
$\hum(r_0;\thetab)=\hat{\hum}_0$.
Unfortunately, there is no available data for the fraction of infected 
mosquitoes. As mosquito abundance is strongly correlated with dengue outbreaks \cite{Jansen2010}, we have chosen a relation 
$\mos(r_0;\thetab)= 3 \hat{\hum}_0$ at the beginning of an epidemic outburst
(other choices gave similar numerical results).

\subsection{Parameter estimation}

To estimate a parameter vector~$\thetab \in \Theta$ in~\eqref{eq:parameter}
that fits with the data provided by the Municipal Secretariat of Public Health 
we apply the curve-fitting approach based on least-square method.
More precisely, we look for an optimal solution to the problem
\begin{equation}
\label{min}
\min \limits_{\thetab \in \Theta} \sum_{r=r_0}^R 
\big( \hum(r;\thetab)-\hat{\hum}_r \big)^2 \eqsepv
R=60~\mbox{days} \eqfinv
\end{equation}
subject to the differential constraint~\eqref{RMmodel-par}. 
Regarding numerics, we have solved this optimization problem with the 
\texttt{lsqcurvefit} routine (MATLAB Optimization Toolbox),
starting with an admissible $\thetab^{(0)} \in \Theta$ 
(the exact value stands in the second column of Table~\ref{tab-A1}). 
The routine generates a sequence 
$\thetab^{(1)}, \thetab^{(2)} \ldots$ that we stop once it is stationary, 
up to numerical precision.
% converges to a viable approximation $\thetab_{\kappa}$ of $\thetab^{*}$ in the sense that $\thetab_{\kappa} \approx \thetab_{\kappa+1}$ within the limits of chosen numerical precision. 
For a better result, we have combined two particular methods 
(Trust-Region-Reflective Least Squares Algorithm \cite{Sorensen1982} and 
Levenberg-Marquardt Algorithm \cite{More1978}) in the implementation of 
the \texttt{lsqcurvefit} MATLAB routine.

  \begin{figure}
  	\begin{center}
	   \includegraphics[width=14cm,height=7cm]{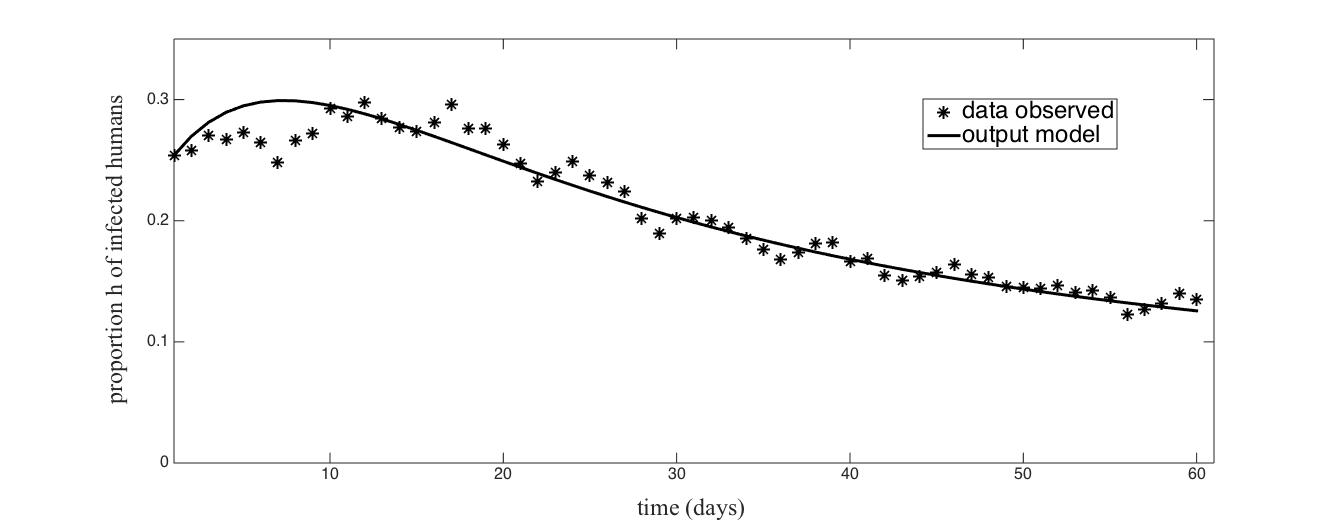}
	\end{center}
	\caption{Fraction of inviduals infected with dengue, obtained by 
adjustment of the Ross-Macdonald model~\eqref{eq:model-Ross-Macdonald} 
(smooth solid curve) versus registered daily prevalence cases 
(star isolated points) during the 2013 dengue outbreak in Cali, Colombia
\label{fig-A1}}
  \end{figure}

The last column of Table~\ref{tab-A1} provides estimated values for the
parameters $\thetab= \big( \alpha, p_{\hum}, p_{\mos}, \xi, \delta \big)$,
and Figure~\ref{fig-A1} displays the curve-fitting results. 
\begin{table}[t]
\begin{center}
\begin{tabular}{|c|c|c|c|c|c|}\hline
Parameter &  Initial & Range & Reference & Estimated & Unit \\ 
&  value &  &  & value & \\ \hline\hline
$\alpha$            &  1     & $[0, 5]$ & \cite{Costero1998}, \cite{Scott2000b}&
0.3600 & $\text{day}^{-1}$ \\ \hline
$p_{\mos}$          &  0.5   & $[0, 1]$ &  & 0.2128 & dimensionless \\ \hline
$p_{\hum}$          &  0.5   & $[0, 1]$ &  & 0.1990 & dimensionless \\ \hline
$\xi$               &  1     & $[1, 5]$ & \cite{Mendez2006}, \cite{Scott2000a} 
& 1.0087 & dimensionless \\ \hline
$\delta$            &  0.035 & $\left[ 0.033, 0.066 \right]$ 
& \cite{Costero1998}, \cite{Scott2000a} & 0.0333 & $\text{day}^{-1}$ \\ \hline
\end{tabular}
\end{center}
\caption{Initial values, admissible ranges, respective source references, 
estimated values of parameters (numerical solution of the
optimization problem~\eqref{min}--\eqref{RMmodel-par}) and units}
\label{tab-A1}
\end{table}
  \begin{subequations}
Therefore, the aggregate parameters~\eqref{eq:transmision-ross} 
%and rate~$\gamma$ of human recovery 
are estimated as
\begin{equation}
A_{\mos}= 0.076608~\text{day}^{-1} \eqsepv A_{\hum}= 0.0722633~\text{day}^{-1} 
%\eqsepv \gamma=0.1~\text{day}^{-1} 
\eqfinp
% 0.3600*0.2128
% 0.3600*0.1990*1.0087
\label{eq:transmision-ross_estimated}
\end{equation}
For the natural mortality rate~$\delta=\lb{\controlm}$ of mosquitoes
(see Table~\ref{tabla:parametros-Ross}), we obtain 
\begin{equation}
  \delta=\lb{\controlm}=  0.0333~\text{day}^{-1} \eqfinp
\label{eq:natural_mortality_rate_estimated}
\end{equation}
For the mosquito mortality maximal rate~$\ub{\controlm}$, we take 
\begin{equation}
 \ub{\controlm}=  0.05~\text{day}^{-1} \eqfinp
\label{eq:maximal_mosquito_mortality_rate_estimated}
\end{equation}
  \end{subequations}

\subsection{Ranges for the uncertain aggregate parameters}
\label{Ranges_for_the_uncertain_aggregate_parameters}

From the ranges for the parameters 
$\big( \alpha, p_{\hum}, p_{\mos}, \xi, \delta \big)$ displayed
in the third column of Table~\ref{tab-A1}, we obtain ranges 
for the aggregate parameters~\eqref{eq:transmision-ross}:
\begin{equation}
A_{\mos}=\alpha p_{\mos} \in [0,5]~\text{day}^{-1}  \eqsepv 
A_{\hum}=\alpha p_{\hum}\xi \in [0,25]~\text{day}^{-1}  \eqfinp
\label{eq:Ranges_for_the_uncertain_aggregate_parameters}
\end{equation}
\bigskip

\paragraph*{Acknowledgments.}

The authors thank the French program PEERS-AIRD 
(\emph{Modèles d'optimisation et de viabilité en écologie et en économie})
and the Colombian Programa Nacional de Ciencias Básicas COLCIENCIAS
(\emph{Modelos y métodos matemáticos para el control y vigilancia del dengue},
código 125956933846)
% Fecha de inicio del proyecto: 18 de febrero de 2014
that offered financial support for missions, together with 
\'Ecole des Ponts ParisTech (France), 
Universit\'e Paris-Est (France), 
Universidad Aut\'onoma de Occidente (Cali, Colombia)
and Universidad del Valle (Cali, Colombia).
We thank the Municipal Secretariat of Public Health 
(Cali, Colombia) for technical support and discussions.
We are indebted to the editor-in-chief and to a reviewer, whose comments
contributed to improve the manuscript.

\newcommand{\noopsort}[1]{} \ifx\undefined\allcaps\def\allcaps#1{#1}\fi


\begin{thebibliography}{10}

\bibitem{AndersonMay1992}
R.~M. Anderson and R.~M. May.
\newblock {\em Infectious Diseases of Humans: Dynamics and Control}.
\newblock Oxford Science Publications. OUP Oxford, 1992.

\bibitem{Aubin1991}
J.~Aubin.
\newblock {\em Viability theory}.
\newblock Systems \& Control: Foundations \& Applications. Birkh\"auser Boston
  Inc., Boston, MA, 1991.

\bibitem{Bene-Doyen:2008}
C.~B\'en\'e and L.~Doyen.
\newblock Contribution values of biodiversity to ecosystem performances: A
  viability perspective.
\newblock {\em Ecological Economics}, 68(1-2):14 -- 23, 2008.

\bibitem{Bene-Doyen-Gabay:2001}
C.~B\'en\'e, L.~Doyen, and D.~Gabay.
\newblock A viability analysis for a bio-economic model.
\newblock {\em Ecological Economics}, 36:385--396, 2001.

\bibitem{Bonneuil-Mullers:1997}
N.~Bonneuil and K.~M\"ullers.
\newblock Viable populations in a prey-predator system.
\newblock {\em Journal of Mathematical Biology}, 35(3):261--293, February 1997.

\bibitem{Bonneuil-Saint-Pierre:2005}
N.~Bonneuil and P.~Saint-Pierre.
\newblock Population viability in three trophic-level food chains.
\newblock {\em Applied Mathematics and Computation}, 169(2):1086 -- 1105, 2005.

\bibitem{BrauerCastillo2006}
F.~Brauer and C.~Castillo-Ch{\'a}vez.
\newblock {\em Mathematical models in population biology and epidemiology},
  volume~40 of {\em Texts in Applied Mathematics}.
\newblock Springer-Verlag, New York, 2001.

\bibitem{Costero1998}
A.~Costero, J.~D. Edman, G.~G. Clark, and T.~W. Scott.
\newblock Life table study of \emph{Aedes aegypti} (diptera: Culicidae) in
  {P}uerto {R}ico fed only human blood plus sugar.
\newblock {\em Journal of Medical Entomology}, 35(5), 1998.

\bibitem{DeLara-Doyen:2008}
M.~{De Lara} and L.~Doyen.
\newblock {\em Sustainable Management of Natural Resources. Mathematical Models
  and Methods}.
\newblock Springer-Verlag, Berlin, 2008.

\bibitem{DeLara-Sepulveda:2016}
M.~{De Lara} and L.~Sepulveda.
\newblock Viable control of an epidemiological model.
\newblock {\em Mathematical Biosciences}, 280:24--37, 2016.

\bibitem{Diekmann-Heesterbeek:2000}
O.~Diekmann and J.~A.~P. Heesterbeek.
\newblock {\em Mathematical Epidemiology of Infectious Diseases}.
\newblock Wiley, Utrecht, Netherland, 2000.

\bibitem{Hethcote:2000}
H.~W. Hethcote.
\newblock The mathematics of infectious diseases.
\newblock {\em SIAM Review}, 42:599--653, 2000.

\bibitem{Jansen2010}
C.~C. Jansen and N.~W. Beebe.
\newblock The dengue vector \emph{Aedes aegypti}: what comes next.
\newblock {\em Microbes and infection}, 12(4):272--279, 2010.

\bibitem{Mendez2006}
F.~M\'endez, M.~Barreto, J.~Arias, G.~Rengifo, J.~Mu\~noz, M.~Burbano, and
  B.~Parra.
\newblock Human and mosquito infections by dengue viruses during and after
  epidemics in a dengue-endemic region of {Colombia}.
\newblock {\em Am J Trop Med Hyg.}, 74(4):678--683, 2006.

\bibitem{More1978}
J.~J. Mor{\'e}.
\newblock The {Levenberg-Marquardt} algorithm: Implementation and theory.
\newblock In G.~A. Watson, editor, {\em Numerical Analysis: Proceedings of the
  Biennial Conference Held at Dundee, June 28--July 1, 1977}, pages 105--116.
  Springer Berlin Heidelberg, Berlin, Heidelberg, 1978.

\bibitem{Regnier-DeLara:2015}
E.~Regnier and M.~{De Lara}.
\newblock Robust viable analysis of a harvested ecosystem model.
\newblock {\em Environmental Modeling \& Assessment}, 20(6):687--698, 2015.

\bibitem{Scott2000b}
T.~W. Scott, P.~H. Amerasinghe, A.~C. Morrison, L.~H. Lorenz, G.~G. Clark,
  D.~Strickman, P.~Kittayapong, and J.~D. Edman.
\newblock Longitudinal studies of \emph{Aedes aegypti} (diptera: Culicidae) in
  {Thailand} and {Puerto Rico}: Blood feeding frequency.
\newblock {\em Journal of Medical Entomology}, 37(1):89, 2000.

\bibitem{Scott2000a}
T.~W. Scott, A.~C. Morrison, L.~H. Lorenz, G.~G. Clark, D.~Strickman,
  P.~Kittayapong, H.~Zhou, and J.~D. Edman.
\newblock Longitudinal studies of \emph{Aedes aegypti} (diptera: Culicidae) in
  {Thailand} and {Puerto Rico}: Population dynamics.
\newblock {\em Journal of Medical Entomology}, 37(1):77, 2000.

\bibitem{Smith_et_al2007}
D.~L. Smith, F.~E. McKenzie, R.~W. Snow, and S.~I. Hay.
\newblock Revisiting the basic reproductive number for malaria and its
  implications for malaria control.
\newblock {\em PLoS Biol}, 5(3):e42, 02 2007.

\bibitem{Sorensen1982}
D.~C. Sorensen.
\newblock Newton's method with a model trust region modification.
\newblock {\em SIAM Journal on Numerical Analysis}, 19(2):409--426, 1982.

\end{thebibliography}
\end{document}